\input amstex
\documentstyle{amsppt}
\magnification=\magstep1
\baselineskip=18pt
\hsize=6truein
\vsize=9truein
\TagsOnLeft

\input epsf
\topmatter

\title On finiteness of Kleinian groups in general dimension
\endtitle

\author Sun-Yung A. Chang \footnote{Research is supported in part by NSF Grant
DMS-0070542 and a Guggenheim Foundation Fellowship.}, Jie Qing
\footnote{Research is supported in part by NSF Grant DMS-9803399
and Slaon Fellowships BR-3818.}  and Paul C. Yang
\footnote{Research is supported in part by NSF Grant DMS-0070526.}
\endauthor

\leftheadtext{On finiteness of Kleinian groups in general
dimension} \rightheadtext{Chang, Qing and Yang}

\abstract In this paper we provide a criteria for geometric
finiteness of Kleinian groups in general dimension. We formulate
the concept of conformal finiteness for Kleinian groups in space
of dimension higher than two, which generalizes the notion of
analytic finiteness in dimension two. Then we extend the argument
in the paper of Bishop and Jones to show that conformal finiteness
implies geometric finiteness unless the set of limit points is of
Hausdorff dimension n. Furthermore we show that, for a given
Kleinian group $\Gamma$,  conformal finiteness is equivalent to
the existence of a metric of finite geometry on the Kleinian
manifold $\Omega(\Gamma)/\Gamma$.
\endabstract

\address Sun-Yung Alice Chang, Department of Mathematics,
Princeton University, \newline
Princeton, NJ 08544 \& Department of Mathematics,
UCLA, Los Angeles, CA 90095. \endaddress
\email chang{\@}math.princeton.edu \endemail
\address Jie Qing, Department of Mathematics, University of California,
Santa Cruz, Santa Cruz,
CA 95064. \endaddress
\email qing{\@}math.ucsc.edu \endemail
\address Paul Yang, School of Mathematics, IAS, Princeton, NJ 08540
\& Department of Mathematics, University of Southern California,
Los Angeles, CA 90089.
\endaddress
\email pyang{\@}math.usc.edu \endemail
\endtopmatter

\document


\vskip .2in \noindent \S0. {\bf Introduction} \vskip .1in

A discrete subgroup $\Gamma$ of the group of conformal
transformations of the unit sphere $S^n$ is called a  Kleinian
group if there is a non-empty domain $\Omega(\Gamma)$ of
discontinuity in $S^n$. The group  $\Gamma$ also acts as a
subgroup of the group of hyperbolic isometries of the unit ball
$B^{n+1}$. Assuming the group has no torsion elements, the
quotient is a hyperbolic manifold $B^{n+1}/\Gamma$, which is
bounded at infinity by a Kleinian manifold $\Omega(\Gamma)/\Gamma$
which is a manifold with a locally conformally flat structure.

There is a useful notion of geometric finiteness for Kleinian
groups that assures nice properties of the geometric quotient.
There are several equivalent forms of this condition. According to
one formulation, a Kleinian group $\Gamma$ is said to be
geometrically finite if the limit set $\Lambda(\Gamma)$ consists
only of conical limit points and cusped limit points. It is a
natural problem to find useful criteria to assure geometric
finiteness. In dimension two, a Kleinian group is said to be
analytically finite if the Riemann surface $\Omega(\Gamma)/
\Gamma$ is of finite type (that is to say a union of a finite
number of closed Riemann surfaces each with finite number of
punctures). Recently Bishop and Jones ([5]) showed that, for an
analytically finite group $\Gamma$, $\Lambda(\Gamma)$ has
Hausdorff dimension less than two if and only if $\Gamma$ is
geometrically finite. An important part of their work is the
construction of an invariant Lipschitz graph which serves to
relate the geometry of the hyperbolic manifold $B^3/\Gamma$ to the
geometry of the Riemann surface $\Omega(\Gamma)/\Gamma$.

Recently we studied locally conformally flat 4-manifolds and
obtained some finiteness for certain class of such manifolds [7]
[8]. In those works [8] the holonomy representation of the
fundamental group of such manifolds as Kleinian group played a key
role in our understanding of the structure of such manifolds. In
dimension higher than two, Kleinian groups have been studied
mostly in conjunction with hyperbolic structure. Our motivation is
to investigate the close relation between the geometry of the
hyperbolic manifolds $B^{n+1}/\Gamma$ and the geometry of the
Kleinian manifold $\Omega(\Gamma)/\Gamma$ for a given Kleinian
group $\Gamma$.

For our purpose, we introduce a notion of conformal finiteness
(see Definition 3.2 in Section 3), which is the natural analogue
of the notion of analytic finiteness to higher dimensions. Then we
extend the theorem of Bishop and Jones to Kleinian groups in
higher dimension.

\proclaim{Theorem 0.1} Suppose that $\Gamma$ is a nonelementary,
conformally finite Kleinian group on $S^n$, then $\Gamma$ is
geometrically finite if and only if the limit set of $\Gamma$ has
Hausdorff dimension strictly smaller than $n$.
\endproclaim

Recall the celebrated finiteness theorem of Ahlfors [2] and Bers
[3] [4], which states that a finitely generated Kleinian group in
dimension two is analytically finite. This finiteness theorem
fails to hold in higher dimensions as pointed out by the examples
of Kapovich [11], and Kapovich and Potyagailo [12]. On the other
hand, a result of Jarvi and Vuorinen [10] shows that in general
dimensions, the limit set of a finitely generated Kleinian group
is uniformly perfect. The latter condition is equivalent, in
dimension two, to the condition of analytic finiteness of the
group. Moreover in [1] Ahlfors showed some weak finiteness for
Kleinian groups in higher dimension: if $\Gamma$ is finitely
generated, then the dimension of the space of certain class of
mixed tensor densities, automorphic under $\Gamma$, is finite (see
also [14] of Hiromi Ohtake). Therefore, it is interesting to
search for the appropriate version of the finiteness result in
higher dimension, particularly for those Kleinian groups with
small limit sets. We take a first step in that direction by
characterizing the conformally finite ends by geometric
conditions. As another application of the Lipschitz graph
construction, we have

\proclaim{Theorem 0.2} Given a Kleinian group $\Gamma$, the
Kleinian manifold $\Omega(\Gamma)/\Gamma$ is of finite geometry
for some metric in the conformal class if and only if $\Gamma$ is
conformally finite.
\endproclaim

By finite geometry here we mean that its curvature and covariant
derivatives of curvature is bounded, and its volume is finite.

We would like to thank Francis Bonahon and Feng Luo for
informative discussions and interest in this work. The second
author would like to thank MSRI for the hospitality. This note is
completed when the second author is visiting MSRI.

\vskip .1in \noindent \S1. {\bf Construction of the Lipschitz
Graph} \vskip .1in

The following construction, by completely elementary means, of the
invariant Lipschitz graph over a domain of discontinuity
$\Omega(\Gamma)$ of a Kleinian group $\Gamma$ is based on the idea
of Bishop and Jones' in [5].  Take a small positive number
$\epsilon_0$ and consider a collection of balls $\{B_\alpha \}$
such that
$$
B_\alpha = B (x_\alpha, d_\alpha) \ \ \text{and} \ \ d_\alpha = \epsilon_0 \cdot
\text{dist}(x_\alpha, L(\Gamma))
\tag 1.1
$$
for each point $x_\alpha \in \Omega(\Gamma)$, where diameters and
distances are all measured
on $S^n$ with the standard metric $g_0$. To construct
an invariant graph we would enlarge the collection to take in
all images of $B_\alpha$ under the group $\Gamma$ and denote the
collection of balls by $B(\Gamma)$. Set
$$
G(\Gamma) = \partial (\bigcup_\beta H_\beta)
\bigcap B^{n+1}
\tag 1.2
$$
where $H_\beta$ is the hyperbolic half space over each
ball $B_\beta$ in $B(\Gamma)$, i.e. the dome whose boundary intersects
perpendicularly at $\partial B_\beta$ with $S^n$. Clearly
$G(\Gamma)$ is a graph over $\Omega(\Gamma)$ in the following
sense
$$
G(\Gamma) = \{ f(x)x: x\in \Omega(\Gamma) \}
\tag 1.3
$$
where $f(x): \Omega(\Gamma) \rightarrow (0, 1)$. In fact
$G(\Gamma)$ is a Lipschitz graph in the sense that
$$
|f(x) - f(y)| \leq M \text{dist}(x,y)
$$
for some $M>0$ and all $x,y \in \Omega(\Gamma)$. Therefore

\proclaim{Proposition 1.1} Given a nonelementary Kleinian group $\Gamma$ and
a small positive number $\epsilon_0$,
the above constructed graph $G(\Gamma)$ is a $\Gamma$-invariant
Lipschitz graph. Moreover
$$
0 < C_1 \leq \frac{1-f(x)}{\text{dist}(x, L(\Gamma))}
\leq C_2
\tag 1.4
$$
for all $x\in \Omega(\Gamma)$, where $C_1, C_2$ only depend on
$\epsilon_0$.
\endproclaim

We need to recall some facts about M\"{o}bius
transformations. For any M\"{o}bius transformation $\gamma$ on
$R^n\cup \{\infty\}, n\geq 3$, we have
$$
|\gamma (x) - \gamma (y)| = |\gamma'(x)|^\frac 12 |\gamma'(y)|^\frac 12 |x-y|
\tag 1.5
$$
(see, for example, (1.3.2) in [13]). In addition if $\gamma$ is
not a composition of just scalings, rotations, or translations,
then
$$
|\gamma'(x)|_e = \frac {1}{\lambda |x - b|^2}
\tag 1.6
$$
for $b\in R^n$ and $\gamma(b) = \infty$ , where $|\cdot|_e$ denotes norm under Euclidean
metric while we will use $|\cdot|_s$ for the norm under the standard metric of the sphere, i.e.
$$
|\gamma'(x)|_s = \frac {1+|x|^2}{1+|\gamma x|^2}|\gamma'(x)|_e.
$$
As a consequence (1.5) and (1.6) we have:

\proclaim{Lemma 1.2} Suppose that $\Gamma$ is a nonelementary Kleinian group, $\Omega(\Gamma)$ is
its domain of discontinuity and $L(\Gamma)= \partial \Omega(\Gamma)$ is its
set of limit points. Then there exists a positive number $C$ such that
$$
\frac 1C \frac {\text{dist}(\gamma(x), L(\Gamma))}{\text{dist}(x, L(\Gamma))}
\leq |\gamma'(x)|_s \leq C \frac {\text{dist}(\gamma(x), L(\Gamma))}
{\text{dist}(x, L(\Gamma))}
\tag 1.7
$$
for all $x\in \Omega(\Gamma)$ and all $\gamma \in \Gamma$.
\endproclaim
\demo{Proof} For any given $x \in \Omega(\Gamma)$ and $\gamma\in \Gamma$. Let
$$
\text{dist}(x, L(\Gamma)) = \text{dist}(x, a), \ \ \text{and} \ \ \text{dist}(\gamma x,
L(\Gamma)) = \text{dist}(\gamma x, b),
$$
for some $a, b \in L(\Gamma)$. Since $L(\Gamma)$ contains more than two points, there
exists a positive number $d_0 >0$ (which is independent of $\{a, b\}$) such that we can
always find a third point $p\in L(\Gamma)$ such that
$$
\text{dist} (p,\ \{a,b\}) \geq d_0
$$
Then
$$
\text{dist}(x, p) \geq \text{dist}(p,a) - \text{dist}(a,x)\geq \text{dist}(p,a) -
\text{dist}(x, p).
$$
Hence,
$$
\text{dist}(x, p) \geq \frac 12 d_0.
$$
Similarly,
$$
\text{dist}(\gamma x, p) \geq \frac 12 d_0.
$$

\vskip 0.2in
\hskip 1.0in\epsfxsize=3.00in \epsfbox{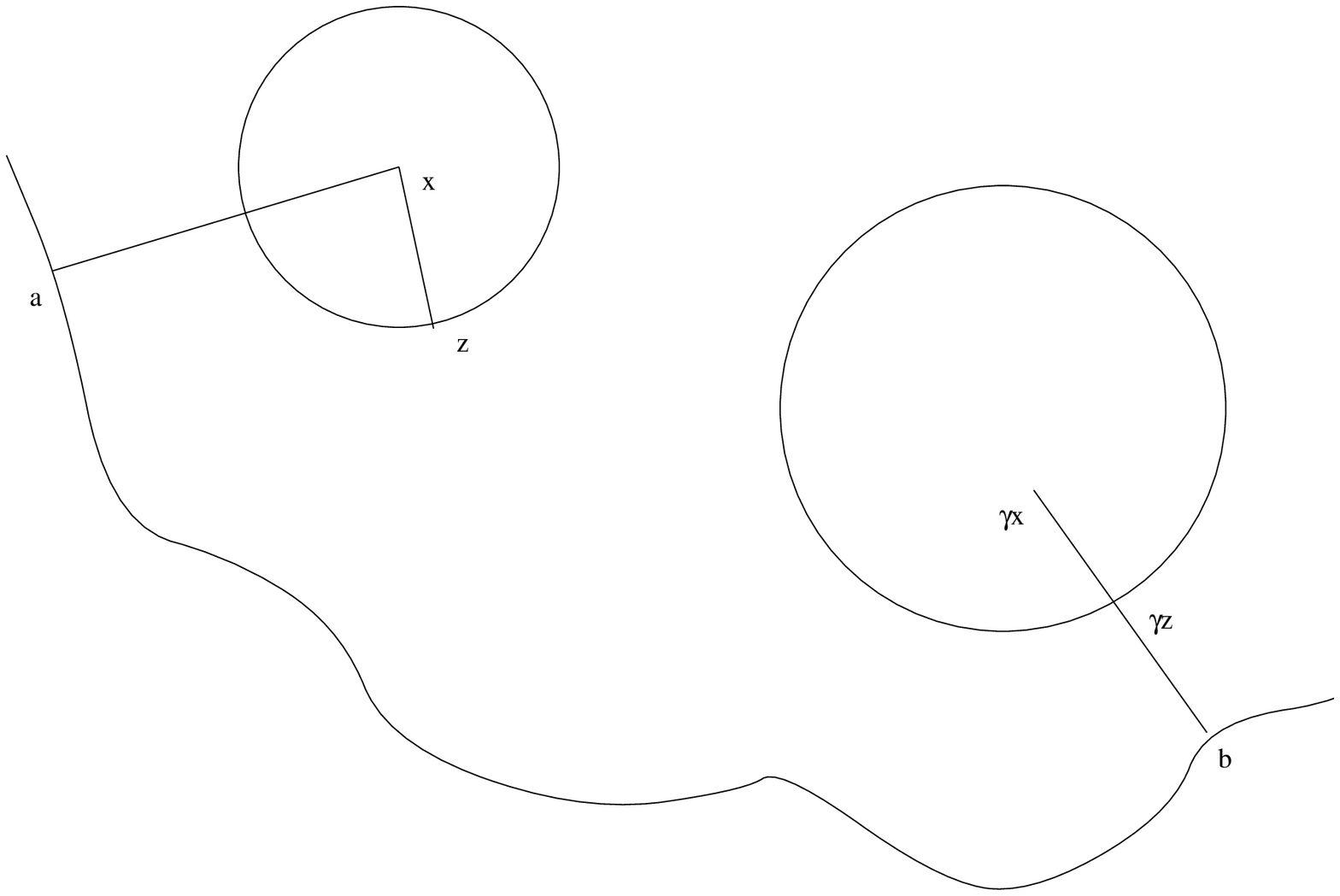}
\vskip 0.2in

Let us fix $p$ as $\{\infty\}$ in $R^n\cup \{\infty\}$. For convenience we continue to use
the same notations for points in $R^n\cup \{\infty\}$. By (1.5) and (1.6) we have
$$
\align
|\gamma x - \gamma z| &= |x-z| |\gamma'(x)|_e^\frac 12 |\gamma'(z)|_e^\frac 12
\\
&= |x-z| \frac 1\lambda \frac 1{|x-c|} \frac 1{|z-c|}
\endalign
$$
for some $c\in L(\Gamma)$ and all $z\in \Omega(\Gamma)$ satisfying $ |x-z| =
\frac 12 |x-a|\}$. Therefore,
$$
|\gamma x - \gamma z| = \frac {|x-z|}{\lambda |x-c|^2} \frac {|x-c|}{|z-c|}
\geq  |x-z||\gamma'(x)|_e \frac {|x-c|}{|x-c| + |x-z|},
$$
Thus
$$
|\gamma x - b| > |\gamma x - \gamma z| \geq \frac 12 |x-a||\gamma'(x)|_e \frac 23.
$$
Due to the choice of $p$ which is some fixed spherical distance away from all four points
$x, \gamma x, a, b$, we have, for some constant $C$ (may depend on $d_0$),
$$
\frac {\text{dist}(\gamma(x), L(\Gamma))}{\text{dist}(x, L(\Gamma))}
\geq \frac 1C |\gamma'(x)|_s.
\tag 1.8
$$
Applying the same argument to $\gamma ^{-1}$ we find
$$
|\gamma'(x)| \geq \frac 1C \frac {\text{dist}(\gamma(x), L(\Gamma))}{\text{dist}(x,
L(\Gamma))}.
\tag 1.9
$$
For those element in the group $\Gamma$ which is a composition of scalings,
rotations, or translations, (1.7) is easily verified.
Therefore the proof is complete.
\enddemo

\noindent {\bf Remark 1.3.} \quad For $L(\Gamma)=\{ \infty \}$ and
$\Gamma$ is simply generated by a translation $\gamma x = x+h$, we
have
$$
|\gamma'(x)|_s = \frac {1+|x|^2}{1+|\gamma x|^2} = \frac {q(\gamma x, \infty)^2}
{q(x, \infty)^2},
$$
where
$$
q(x,y) = \frac {2|x-y|}{\sqrt{1+|x|^2}\sqrt{1+|y|^2}}, \ \ \ \text{and} \ \ \ \
q(x, \infty) = \frac 2{\sqrt{1+|x|^2}}
$$
which is called the chordal metric and is equivalent to the
spherical distance $d(x,y)$. For $L(\Gamma)= \{0, \infty\}$ and
$\Gamma$ is simply generated by an inversion $\gamma x = \frac
x{|x|^2}$, we also have
$$
|\gamma' (x)|_s = \frac {q(\gamma x, \infty)^2}{q(x, 0)^2}.
$$
\vskip0.1in 

\demo{Proof of Proposition 1.1}
 It is clear from the construction that
the graph $G(\Gamma)$ is $\Gamma$-invariant. We need to verify that it is a Lipschitz graph.
This follows easily as long as for all balls in $B(\Gamma)$ we have
$$
C_0 \leq \frac {\text{diam}(D_\beta)}{\text{dist}(D_\beta, L(\Gamma))} \leq \frac 12,
\tag 1.10
$$
for some positive constant $C_0$. But (1.10) is a consequence of (1.7) in Lemma 1.2 for
some appropriate choice of $\epsilon_0$. So the proof is complete.
\enddemo

Given a Kleinian group $\Gamma$, let $C(L(\Gamma)) \subset
B^{n+1}$ denote the convex hull of the set of limit points
$L(\Gamma) \subset S^n$ with respect to the hyperbolic metric on
$B^{n+1}$. Then $C(B^{n+1}/\Gamma) = C(L(\Gamma))/\Gamma$ is
called the convex core of the hyperbolic manifold
$B^{n+1}/\Gamma$. From the construction of the graph $G(\Gamma)$
the following is clear.

\proclaim{Lemma 1.4} For any nonelementary Kleinian group, the
above constructed $\Gamma$-invariant graph $G(\Gamma)$ separates
the convex hull $C(L(\Gamma))$ from $\Omega(\Gamma)$.
\endproclaim

\vskip 0.2in \noindent \S2 {\bf Compact Case} \vskip .1in The
Lipschitz graph constructed in the previous section has an induced
metric from the hyperbolic metric. We would like to compare it
with a suitable conformal metric on the Kleinian quotient
$\Omega(\Gamma)/\Gamma$. In the case when the Kleinian quotient is
compact, this is a relatively simple matter.

\proclaim{Theorem 2.1} Suppose that $\Gamma$ is a nonelementary
Kleinian group, and that the Kleinian manifold
$\Omega(\Gamma)/\Gamma$ is compact. Then, for any metric in the
conformal class on $\Omega(\Gamma)/\Gamma$, we have a complete
$\Gamma$-invariant metric $e^{2u}g_0$ on $\Omega(\Gamma)$ and
$$
\frac 1K \frac 1{\text{dist}(x, L(\Gamma))} \leq e^{u(x)} \leq K \frac 1{\text{dist}(x,
L(\Gamma))},
\tag 2.1
$$
for all $x\in \Omega(\Gamma)$ and some positive number $K$.
\endproclaim

\demo{Proof} This basically is a consequence of Lemma 1.2. Due to the invariance of
the metric $e^{2u}g_0$, we have
$$
e^{u(x)} = e^{u(\gamma x)} |\gamma'(x)|_s.
\tag 2.2
$$
Now, fix a fundamental region $D$, whose closure $\bar D$
is compact in $\Omega(\Gamma)$, and for any
$x\in \Omega(\Gamma)$, there exists $\gamma\in \Gamma$ such that
$\gamma^{-1} x =y \in \bar D$. Then
$$
e^{u(y)} = e^{u(x)} |\gamma' (y)|_s
$$
where, by Lemma 1.2,
$$
\frac 1C \frac {\text{dist}(x, L(\Gamma))}{\text{dist}(y, L(\Gamma))}
\leq |\gamma' (y)|_s \leq C \frac {\text{dist}(x, L(\Gamma))}
{\text{dist}(y, L(\Gamma))}.
$$
Thus
$$
\frac 1C \frac {e^{u(y)} \text{dist}(y, L(\Gamma))}{\text{dist}(x, L(\Gamma))}
\leq e^{u(x)} \leq C\frac {e^{u(y)} \text{dist}(y, L(\Gamma))}{\text{dist}(x, L(\Gamma))},
\tag 2.3
$$
that is
$$
\frac 1K \frac 1{\text{dist}(x, L(\Gamma))} \leq e^{u(x)}
\leq K \frac 1{\text{dist}(x, L(\Gamma))}
$$
for some positive constant $K$.
\enddemo

\noindent {\bf Remark 2.2.} \quad In particular, we may consider
the Yamabe metric on $\Omega(\Gamma)/ \Gamma$ in Theorem 2.1.
Therefore, under the assumptions of Theorem 2.1, there is a
complete $\Gamma$-invariant metric on $\Omega(\Gamma)$ with
constant scalar curvature and satisfying (2.1). \vskip 0.1in

We point out that the natural bounds (2.1) on the invariant metric on
$\Omega(\Gamma)$ is the key to relate the hyperbolic
geometry inside $B/\Gamma$ to the conformal geometry at infinity $\Omega(\Gamma)/\Gamma$
through the constructed $\Gamma$-invariant Lipschitz graph $G(\Gamma)$. Namely,

\proclaim{Proposition 2.3} Suppose that $\Gamma$ is a nonelementary Kleinian group and that
the Kleinian manifold $\Omega(\Gamma)/\Gamma$ is compact. Then the map
$$
F(x) = f(x) \, x : \Omega(\Gamma) \rightarrow G(\Gamma)
$$
is a $\Gamma$-invariant bi-Lipschitz map with respect to the induced hyperbolic metric on
the graph $G(\Gamma)$ and any metric on $\Omega(\Gamma)$ which is induced by a metric on
$\Omega(\Gamma)/\Gamma$ in the conformal class.
\endproclaim

\vskip 0.2in
\noindent
\S3 {\bf Geometric finiteness}
\vskip .1in

In this section we discuss notions of finiteness for Kleinian
groups. We begin with recalling notions of geometric finiteness
for Kleinian groups in the study of hyperbolic manifolds. Then we
give a definition of conformal  finiteness as a generalization of
analytic finiteness for Kleinian groups as discrete subgroups of
conformal transformations on $S^2$. Then we will give some useful
metrics on the Kleinian manifold $\Omega(\Gamma)/\Gamma$ for
conformally finite group $\Gamma$.

A good reference for the discussion of geometric finiteness for
Kleinian groups is the paper of Bowditch [6]. We also refer
readers to the book of Ratcliffe [15] and Alhfors' lecture notes
[1] for all basics about Kleinian groups.

\vskip 0.1in \noindent {\bf Definition 3.1.} \quad A Kleinian
group $\Gamma$, as a discrete subgroup of the group of hyperbolic
isometries, is geometrically finite if its limit set $L(\Gamma)$
consists entirely of conical limit points and cusped limit points.
\vskip 0.1in

An equivalent formulation says  $\Gamma$ is geometrically finite if and only if the
thick part of the convex
core $C(B^{n+1}/\Gamma)$ is compact. Thus, if $\Gamma$ is geometrically infinite,
 there must be a sequence of points $\{p_i\} \in C(B^{n+1}/\Gamma)$ for which
the injective radius of $B^{n+1}/\Gamma$ at $p_i$ is bounded from
below and $p_i$ tends to infinity in the convex core
$C(B^{n+1}/\Gamma)$. This fact will be  used later. Another
equivalent definition of geometric finiteness says $\Gamma$ is
geometrically finite if and only if $(B^{n+1}\cup
\Omega(\Gamma))/\Gamma$ may be considered as a union of a compact
set and a finite number of disjoint standard cusp ends. It follows
that $\Omega(\Gamma)/\Gamma$ is the union of a compact set and a
finite number of disjoint ends which we will call standard
conformal cusp ends. Of course, a standard conformal cusp end
$C_m$ is the ideal boundary of the standard hyperbolic cusp end.
For a discrete subgroup $\Gamma_\infty$ of the group of Euclidean
isometries of $R^n$, let $R^{n-m}$ be the maximal invariant
subspace so that $R^{n-m}/\Gamma_\infty$ is compact. Suppose that
$N (R^{n-m}, \epsilon)$ is an $\epsilon$-neighborhood of $R^{n-m}$
in $R^n$. Then $N (R^{n-m}, \epsilon)$ is also invariant under
$\Gamma_\infty$ and a standard conformal cusp end is of the form
$$
(R^n \setminus N (R^{n-m}, \epsilon))/\Gamma_\infty.
\tag 3.1
$$
Therefore, a standard conformal cusp end is conformal to
$(R^m\setminus B_R(0)) \times K$ where $K$ is a compact locally
flat manifold of dimension $n-m$. Now we are ready to give the
following definition.

\vskip 0.1in \noindent {\bf Definition 3.2.} \quad Suppose that
$\Gamma$ is a Kleinian group. Then we say $\Gamma$ is conformally
finite if $\Omega(\Gamma)/\Gamma$ is a disjoint union of a compact
set and a finite number of standard conformal cusp ends. \vskip
0.1in

By definition, geometric finiteness implies conformal
finiteness. In this terminology, our goal is to
investigate when conformal finiteness implies geometric finiteness.
It is
clear that the notion of conformal finiteness is a higher dimension analogue of
the analytic finiteness.

In the case that $\Gamma$ is an analytic finite group acting on $S^2$, the uniformization
theorem yields a  hyperbolic
metric on $\Omega(\Gamma)/\Gamma$ and the Lemma of Schwarz and the Koebe
distortion theorem show that the hyperbolic metric satisfies the bounds (2.1).
In higher
dimension there are several possible canonical metrics available, but the
natural bounds (2.1) becomes an issue.
When $\Gamma$ is conformally finite,  we describe a metric
on the Kleinian manifold $\Omega(\Gamma)/\Gamma$ by  explicitly writing metrics on
each conformal cusp end. Take $(x,y)\in R^m\times R^{n-m}$ and let
$$
g_h = \frac 1{|x|^2}(|dx|^2 + |dy|^2).
$$
The set $\{(x,y)|x \neq 0\}$ is conformally the standard
  $H^{n-m+1}\times S^{m-1}$,
it is the holonomy cover of our typical conformally finite end.

\proclaim{Lemma 3.3} $g_h$ induces a complete metric in the standard conformal
class of $(R^m\setminus \{0\})\times K$ with constant scalar curvature
$\frac {n-2}4(2m-n-2)$. More importantly, if we write $g_h = e^{2u}g_0$ where $g_0$ is
the standard metric on $S^n$, then
$$
\frac 1C \frac 1{\text{dist}(p, \infty)} \leq
e^{u(p)} \leq C \frac 1{\text{dist}(p, \infty)}
\tag 3.2
$$
for some constant $C>0$ and all $p$ in a fundamental domain $R^m\times K'$ where
$\bar {K'}/\Gamma_\infty = K$. A standard conformal cusp end
$((R^m\setminus B_R(0))\times K, g_h)$ has a finite volume.
\endproclaim

\demo{Proof}
Let us verify (3.2), we
again use the chordal distance instead of spherical distance. Recall
$$
q((x,y), \infty) = \frac 2{\sqrt{1+|x|^2+|y|^2}}.
$$
Therefore, for $p=(x,y)$,
$$
q(p, \infty)e^{u(p)} = \frac {\sqrt{1+|x|^2+|y|^2}}{|x|}.
\tag 3.3
$$
Thus, when restrict $p$ in a fundamental domain, for instance, $R^m\times K'$, we
have $|y|$ bounded and
$$
\frac 1C' \leq q(p, \infty)e^{u(p)} \leq C'
\tag 3.4
$$
for some constant $C'$ depending on the size of the fundamental domain of $\Gamma_\infty$,
which implies (3.2). To compute the volume, we have
$$
\aligned
\text{vol} (R^m\setminus B_R(0))\times K)
& = \int_{R^m\setminus B_R(0)}\int_{K'} |x|^{-n} dydx \\
& = \text{vol}(K)\int_{R^m\setminus B_R(0)} |x|^{-n}dx \\
& = \text{vol}(K)\int_R^\infty \frac 1{t^{n-m+1}}dt \\
& = \text{vol}(K) \frac 1{n-m} \frac 1{R^{n-m}}.
\endaligned
\tag 3.5
$$
So the proof is complete.
\enddemo

Now let us fix a conformal metric on $\Omega(\Gamma)/\Gamma$ which agrees with
the $g_h$ on each conformal cusp end and arbitrary on the compact part. Let us denote it by
$g_\Gamma$ (this is not intended to signify $g_{\Gamma}$ is an any way canonical).

\proclaim{Lemma 3.4} Suppose that $\Gamma$ is nonelementary,
conformally finite Kleinian group, and that
$g_\Gamma$ is a metric constructed as the above. Then the metric $e^{2u}g_0$
on $\Omega(\Gamma)$ lifted from $g_\Gamma$ satisfies
$$
\frac 1C \frac 1{\text{dist}(x, L(\Gamma))} \leq e^{u(x)}
\leq C \frac 1{\text{dist}(x, L(\Gamma))}
\tag 3.6
$$
for some constant $C>0$ and all $x \in \Omega(\Gamma)$.
\endproclaim

\demo{Proof} In light of (2.3) in Section 2, we only need to verify (3.6) for
all $x$ in a fundamental region. In a fundamental region for a conformally
finite $\Gamma$, one only needs to verify (3.6) for all $x$ in the part corresponding to
each conformal cusp end, which is supported by (3.2) in Lemma 3.3. Thus the proof is complete.
\enddemo

To relate the geometry of the hyperbolic manifold $B^{n+1}/\Gamma$
and the Kleinian manifold $\Omega(\Gamma)/\Gamma$, we find

\proclaim{Proposition 3.5} Suppose that $\Gamma$ is nonelementary,
conformally finite, and that
$g_\Gamma$ is a metric constructed as the above. Then the map
$$
F(x) = f(x)x: \Omega(\Gamma) \rightarrow G(\Gamma)
$$
is a $\Gamma$-invariant bi-Lipschitz map with respect to the induced hyperbolic metric
on the graph $G(\Gamma)$ and the above metric $g_\Gamma$ on $\Omega(\Gamma)$. Moreover
the hypersurface $G(\Gamma)/\Gamma$ in $B^{n+1}/\Gamma$ has finite volume.
\endproclaim

\vskip .2in \noindent \S4 {\bf Hyperbolically Harmonic Functions}
\vskip .1in

This section is concerned with harmonic functions on the
hyperbolic manifolds $B^{n+1}/\Gamma$ for a given Kleinian group
$\Gamma$. Good references are Chapter V in Ahlfors' lecture notes
[1] and Chapter V in Nicholls' book [13]. We begin with Green's
function on $B^{n+1}/\Gamma$. In this article we are only
concerned with the Kleinian group of second kind, which simply
means $L(\Gamma)\neq S^n$. According to Lemma 2 and Theorem 1 in
Chapter VI in Ahlfors' Lecture notes [1], $B^{n+1}/\Gamma$ always
possesses a unique minimal positive Green's function, which is of
the form
$$
G(x,y) = \sum_{\gamma \in \Gamma} g(x, \gamma y)
\tag 4.1
$$
where
$g(x,y) $ is the Green's function on the hyperbolic $B^{n+1}$ and
$$
g(x,y) = g(0, |T_x y|) = \int_{|T_x y|}^1 \frac {(1-t^2)^{n-1}}{t^n}dt
\tag 4.2
$$
(see, Chapter V in Alhfors' lecture notes [1]).  On the other
hand, the Green's function on $B^{n+1}/\Gamma$ is also obtained by
integrating the heat kernel:
$$
G(x,y) = \int_0^\infty H(x,y,t) dt.
\tag 4.3
$$
Therefore, based on bounds for the heat kernel given in Davies
[9], we have the following upper bound for the Green's function.

\proclaim{Lemma 4.1} Suppose that the first eigenvalue $\lambda_0$ of $B^{n+1}/\Gamma$
is positive. Then, for $0< \delta < \frac 12\lambda_0$, we have
$$
0 < G(x,y) < C \text{vol}(B_1(x))^{-\frac 12}\text{vol}(B_1(y))^{-\frac 12}
e^{- \sqrt{\frac {4(\lambda_0-2\delta)}{4+\delta}}\rho(x,y)},
\tag 4.4
$$
where $\rho(x,y)$ is the hyperbolic distance, for some constant $C>0$ and $\rho(x,y)
> 32\sqrt{\lambda_0}$.
\endproclaim

\demo{Proof} Recall that
$$
0 \leq H(x,y,t) \leq C \text{vol}(B_t(x))^{-\frac 12}\text{vol}(B_t(y))^{-\frac 12}
e^{-\frac{\rho(x,y)^2}{(4+\delta)t}}
\tag 4.5
$$
for all $0<t<1$, and
$$
0\leq H(x,y,t) \leq C \text{vol}(B_1(x))^{-\frac 12}\text{vol}(B_1(y))^{-\frac 12}
e^{-(\lambda_0-\delta)t} e^{-\frac {\rho(x,y)^2}{(4+\delta)t}}
\tag 4.6
$$
for all $1 \leq t <\infty$. First, by the nonincreasing property of the function
$$
t^{-(n+1)} e^{-nt} \text{vol}(B_t(x))
$$
proved in Proposition 4.3 in Chapter 1 of [16] we have
$$
H(x,y,t) \leq C \text{vol}(B_1(x))^{-\frac 12}\text{vol}(B_1(y))^{-\frac 12}
t^{-(n+1)}e^{-\frac{\rho(x,y)^2}{(4+\delta)t}}
\tag 4.7
$$
for all $0<t<1$. Therefore
$$
G(x,y)  =  \int_0^1 H(x,y,t)dt + \int_1^\infty H(x,y,t)dt
$$
$$
\aligned
\leq C \text{vol}(B_1(x))^{-\frac 12} & \text{vol}(B_1(y))^{-\frac 12}\{
\int_0^1 t^{-(n+1)} e^{- \frac {\lambda_0} t}dt + \int_1^\infty e^{-\delta t}dt\}\\
& \cdot  e^{-\sqrt{\frac {4(\lambda_0-2\delta)}{4+\delta}}\rho(x,y)}
\endaligned
$$
where we use the fact that
$$
(\lambda_0 - 2\delta)t + \frac {\rho(x,y)^2}{(4+\delta)t} \geq \sqrt{\frac {4(\lambda_0
-2\delta)}{4+\delta}} \rho(x,y).
$$
So the proof is complete.
\enddemo

Now, let us consider the harmonic function on the hyperbolic
$B^{n+1}$ with Dirichlet boundary condition on $S^n$. Given a
$f\in L^1(S^n)$, according to Ahlfors (Chapter V of Ahlfors'
lecture notes [1]), we have a harmonic function $u$ on $B^{n+1}$
as
$$
u(x) = \frac 1{\text{vol}(S^n)} \int_{S^n} k(x,y)^n f(y)d\omega (y),
\tag 4.8
$$
where $d\omega$ is the standard volume element for the sphere $S^n$.
So, if we consider $\chi_{\Omega(\Gamma)}$: the characteristic function of the domain
of discontinuity of $\Gamma$, then, its harmonic extension
$$
u_\Gamma (x) = \frac 1{\text{vol}(S^n)} \int_{S^n} k(x,y)^n \chi_{\Omega(\Gamma)}(y)
d\omega(y)
\tag 4.9
$$
is $\Gamma$-invariant, therefore descends to a harmonic function on the hyperbolic
manifold $B^{n+1}/\Gamma$. We denote by $\omega_\Gamma (x)$ the descended harmonic function
on $B^{n+1}/\Gamma$. Recall, for a nonelementary Kleinian group $\Gamma$, we have constructed
a $\Gamma$-invariant Lipschitz graph $G(\Gamma)$ over $\Omega(\Gamma)$. Thus we have
a hypersurface $S_\Gamma = G(\Gamma)/\Gamma$ in the hyperbolic manifold $B^{n+1}/\Gamma$
separating the convex core $C(B^{n+1}/\Gamma)$ from the ideal boundary $\Omega(\Gamma)/\Gamma$.
We have the following representation formula:

\proclaim{Lemma 4.2} Suppose that $\Gamma$ is nonelementary Kleinian group. Then
$$
\omega_\Gamma (x) = \frac 1{2^{n-1}
\text{vol}(S^n)}\int_{S_\Gamma} ( - \frac
{\partial G}{\partial n} (x,y)) d\sigma (y)
\tag 4.10
$$
where $\frac {\partial}{\partial n}$ is the hyperbolic normal derivative of the hypersurface
$S_\Gamma$ in $B^{n+1}/\Gamma$, and $d\sigma$ is the induced one from $B^{n+1}/\Gamma$.
\endproclaim

\demo{Proof} The proof given by Bishop and Jones in [5] works even
in higher dimension with little modifications. But for the
convenience of the reader, we present the proof here. We start
with $x=0$, namely,
$$
\omega_\Gamma (0) = u_\Gamma (0) = \frac 1{\text{vol}(S^n)}
\int_{\Omega(\Gamma)}d\omega = \frac {\text{vol}(\Omega(\Gamma))}{\text{vol}(S^n)}.
\tag 4.11
$$
Recall that
$$
- \frac{\partial g}{\partial n}d\sigma |_{\partial B_r(0)} =
\frac {2^{n-1}}{r^n}d\omega|_{B_r(0)}.
\tag 4.12
$$
Let $\Omega_r$ be the part of $\partial B_r(0)$ which is between
$G(\Gamma)$ and $\Omega(\Gamma)$ and $G_r = G(\Gamma)\cap B_r$. Clearly
$$
\aligned
u_\Gamma (0) & = \frac 1{\text{vol}(S^n)} \lim_{r \rightarrow 1}
\frac {\text{vol}(\Omega_r)}{r^n} \\
& = \frac 1{2^{n-1}\text{vol}(S^n)} \lim_{r \rightarrow 1}
\int_{\Omega_r} (-\frac{\partial g}{\partial n})d\sigma \\
& = \frac 1{2^{n-1}\text{vol}(S^n)} \lim_{r \rightarrow 1}
\int_{G_r} (-\frac {\partial g}{\partial n})d\sigma
\endaligned
$$
by the fact that $g$ is harmonic in the region bounded by $G_r$ and
$\Omega_r$. Thus
$$
u_\Gamma (0) = \frac 1{2^{n-1}\text{vol}(S^n)}
\int_{G(\Gamma)} (-\frac {\partial g}{\partial n})d\sigma.
\tag 4.13
$$
Notice here that we have used the fact
that the constructed $\Gamma$-invariant graph $G(\Gamma)$ is
Lipschitz, i.e.
$$
\int_{G(\Gamma)}|\frac {\partial g}{\partial n}|d\sigma < \infty.
\tag 4.14
$$
Then, by dominated convergence theorem and (4.1), we have, if let
$S$ be any fundamental region for $G(\Gamma)$,
$$
\aligned
& \int_{G(\Gamma)} (-\frac {\partial g}{\partial n})d\sigma
= \sum_{\gamma \in \Gamma} \int_{\gamma S}
(-\frac {\partial g}{\partial n})d\sigma \\
= & \sum_{\gamma \in \Gamma} \int_S (-\frac {\partial
g(\gamma 0,\gamma y)}{\partial n})d\sigma
= \int_S (-\frac {\partial \sum_{\gamma \in \Gamma}
g(\gamma 0,\gamma y)}{\partial n})d\sigma \\
= & \int_{S_\Gamma} (-\frac{\partial G(0,y)}{\partial n})d\sigma.
\endaligned
\tag 4.15
$$
This proves the lemma for $x=0$. For a general point $x \in S_\Gamma$, take
a conformal transformation $T$ of $S^n$ such that $T0 =x$. Then
$$
\omega_\Gamma (x) = u_\Gamma (T0)
= \frac 1{\text{vol}(S^n)}\int_{S^n} \chi_{T\Omega(\Gamma)} d\omega
= \frac {\text{vol}(T\Omega(\Gamma))}{\text{vol}(S^n)}.
\tag 4.16
$$
Therefore, similarly, we can verify (4.10) for all $x\in S_\Gamma$.
\enddemo

\noindent {\bf Remark 4.3.} \quad The constant $\frac
1{2^{n-1}\text{vol}(S^n)}$ in the formula (4.10) depends on the
choice of the Green function $g(x,y)$.

\vskip 0.2in \noindent \S5 {\bf Proof of Theorem 0.1} \vskip.1in

We begin with the thick-thin decomposition for a hyperbolic space
$B^{n+1}/\Gamma$. Good references are Section 12.5 in Ratcliff's
book [15] and Section 3.3 in Bowditch's paper [6]. We recall that,
by Margulis lemma, there is a dimensional constant $c_n>0$ such
that, for any $\epsilon < c_n$,
$$
V(\Gamma, \epsilon) = \{x \in B^{n+1}: d_H(x, \gamma x) < \epsilon, \ \ \text{for some
$\gamma\in \Gamma$}\}
\tag 5.1
$$
is a disjoint union of connected components
$$
V(\Gamma_a, \epsilon) = \{x\in B^{n+1}: d_H(x, \gamma x) < \epsilon, \ \ \text{for some
$\gamma\in \Gamma_a$}\}
\tag 5.2
$$
where $\Gamma_a$ is either a maximum parabolic elementary subgroup or a maximum
hyperbolic elementary subgroup of $\Gamma$.
Each connected component $V(\Gamma_a, \epsilon)$ for $\epsilon < c_n$ is called a Margulis
region. Notice that
$$
\gamma V(\Gamma_a, \epsilon) \bigcap V(\Gamma_a, \epsilon) = \emptyset
\tag 5.3
$$
for all $\gamma \in \Gamma \setminus \Gamma_a$. Because, otherwise, for some $x\in V(\Gamma_a,
\epsilon)$ and some $\gamma_a \in \Gamma_a$,
$$
d_H(\gamma x, \gamma_a \gamma x) < \epsilon.
$$
Then
$$
d_H(x, \gamma^{-1}\gamma_a\gamma x)< \epsilon,
$$
which implies, by Margulis lemma, $\gamma^{-1}\gamma_a\gamma \in \Gamma_a$, i.e.
$\gamma_a\gamma a = \gamma a$. So $\gamma \in \Gamma_a$. This proves (5.3). Therefore
the thin part $V(\Gamma, \epsilon)/\Gamma$ is disjoint union of connected components
where each component has the form $V(\Gamma_a, \epsilon)/\Gamma_a$. A component
$V(\Gamma_a, \epsilon)/\Gamma_a$ is called a Margulis cusp if $\Gamma_a$ is parabolic
otherwise called a Margulis tube. Since Margulis cusps and Margulis tubes represent
the thin part of the hyperbolic manifold $B^{n+1}/\Gamma$ it is difficult to relate them to
the Kleinian manifold $\Omega(\Gamma)/\Gamma$ directly, in contrast to the
standard cusped region for a cusped limit point. The Lipschitz
hypersurface $S_\Gamma$ is designed to make this comparison possible.

\proclaim{Lemma 5.1} Suppose that $\Gamma$ is nonelementary, conformally
finite Kleinian group, then
$$
\int_{S_\Gamma} \text{vol}(B_1(x))^{-\frac 12}d\sigma (x) < \infty.
\tag 5.4
$$
\endproclaim

\demo{Proof}
Let $\epsilon$ be chosen to be small than the Margulis constant $c_n$, then
we decompose the hypersurface into the Margulis region
$S_\Gamma \cap V(\Gamma, \epsilon)/\Gamma$ and its complement $S'$ which
is compact.
Then clearly
$$
\int_{S_\Gamma} \text{vol}(B_1(x))^{-\frac 12}d\sigma (x)
\leq C + \int_{S_\Gamma \cap V(\Gamma, \epsilon)/\Gamma}
\text{vol}(B_1(x))^{-\frac 12}d\sigma (x).
\tag 5.5
$$
In fact, for any Margulis tube $M_a = V(\Gamma_a,
\epsilon)/\Gamma_a$, it is easily seen that $M_a\cap S_\Gamma$ is
compact. Otherwise, the fixed point $a$ of a hyperbolic subgroup
$\Gamma_a$ would be on the boundary of a fundamental region for
$\Gamma$ on $S^n$, which is impossible. Now suppose that
$\Omega(\Gamma)/\Gamma = M_c \cup (\cup_k C_k)$ where $M_c$ is
compact and $\{C_k\}$ are finite number of conformal cusp ends.
Suppose $a_k$ is the parabolic fixed point associated with the
conformal cusp end $C_k$. Then the only Margulis cusps that has
noncompact intersection with $S_\Gamma$ are those which is
associated with the parabolic fixed point $a_k$ and its stabilizer
$\Gamma_{a_k}$. Thus
$$
\int_{S_\Gamma} \text{vol}(B_1(x))^{-\frac 12}d\sigma (x)
\leq C + \sum_k \int_{S_\Gamma \cap (V(\Gamma_{a_k}, \epsilon)/\Gamma_{a_k})}
\text{vol}(B_1(x))^{-\frac 12}d\sigma (x)
\tag 5.6
$$
$$
\leq C + \sum_k \int_{R^{m_k}\setminus B_1} \int_{K_k} \text{vol}(B_1(F^{-1}(x, y))^{-\frac 12}
|x|^{-n}dy dx
\tag 5.7
$$
where $F: G(\Gamma) \rightarrow \Omega(\Gamma)$ and
$C_k$ is conformal to $K_k\times R^{m_k}$ with the metric as $g_h$ as given in Section 3.
So
$$
\aligned
& \leq C + C \sum_k \text{vol}(K_k)
\int_1^\infty |x|^{\frac {n-m_k}2} |x|^{-n} |x|^{m_k -1}d|x| \\
& \leq C + C \sum_k \int_1^\infty |x|^{-\frac {n-m_k}2 -1}d|x| \\
& < \infty
\endaligned
\tag 5.8
$$
since $m_k \leq n-1$. Therefore the proof is complete.
\enddemo

The following lemma is adopted from Bishop and Jones' paper (see,
Lemma 3.6 in [5]). Their proof applies to higher dimension with
little modifications.

\proclaim{Lemma 5.2} Suppose that $\Gamma$ is nonelementary, conformally
finite Kleinian group and suppose that $\Gamma$ is geometrically infinite. Then, there exist
$\epsilon >0$ and
a sequence of points $x_n\in C(B^{n+1}/\Gamma)$
such that $d_H (x_n, S_\Gamma) \rightarrow \infty$ and the injectivity radius
$\text{inj}(x_n) > \epsilon$ for all $n$.
\endproclaim

\demo{Proof} When $\Gamma$ is conformally finite, the
noncompact part of the surface $S_\Gamma$ has to be inside the disjoint union of
finite number of Margulis cusps, where the injectivity radius decays exponentially
as the point moves towards infinity. Therefore the sequence of points tending to infinity
in the thick part of the hyperbolic manifold must move away from the hypersurface.
\enddemo


With all preparations in place, we are now ready to state and prove the main theorem.

\proclaim{Theorem 5.4} Suppose that $\Gamma$ is nonelementary,
conformally finite Kleinian group, then it is geometrically finite
unless the Hausdorff dimension of its limit point set on $S^n$ is
$n$.
\endproclaim

\demo{Proof} We will follow the argument in Bishop and Jones'
paper [5]. We consider the harmonic function $\omega_\Gamma$
discussed in Lemma 4.2. Then
$$
\omega_\Gamma (x) = C \int_{S_\Gamma} \frac {\partial }{\partial n} G(x,y) d\sigma (y).
\tag 5.9
$$
Since $G(x,y)$ is a harmonic function of $y$ away from $x$.
Therefore, by the gradient estimate (see, for example, Corollary
3.2 in [16]), we have
$$
|\frac {\partial }{\partial n} G(x,y)| \leq C G(x,y).
\tag 5.10
$$
Before we apply the estimate for the Green function in Lemma 4.1,
we notice that, first, Bishop and Jones proves in [5] that
$$
\delta (\Gamma) = \text{dim}(L_c(\Gamma))
\tag 5.11
$$
for all nonelementary Kleinian group $\Gamma$, where
$\delta(\Gamma)$ is the Poincar\'{e} exponent and $L_c(\Gamma)$ is
the set of all conical limit points (they only state (5.11) in
2-dimension, but as they pointed out their argument proves (5.11)
in higher dimension too.); second, Sullivan  [17] generalized
Elstrodt-Patterson theorem in higher dimension as
$$
\lambda_0 (B^{n+1}/\Gamma) = \left\{\aligned
(\frac n2)^2 \ \ \ \ \ & \text{if $\delta(\Gamma) \leq \frac n2$}. \\
\delta(\Gamma)(n - \delta(\Gamma)) \ \ \ \ & \text{if $\delta(\Gamma) \geq \frac n2$}.
\endaligned\right.
\tag 5.12
$$
Therefore $\lambda_0 >0$. By (4.4) and (5.10) we arrive at
$$
\omega_\Gamma (x) \leq C e^{-\sqrt{\frac {4(\lambda_0 - 2\delta)}{4+\delta}}d_H(x,y)}
\int_{S_\Gamma} \text{vol}(B_1(x))^{-\frac 12}\text{vol}(B_1(y))^{-\frac 12}d\sigma(y)
\tag 5.13
$$
Now, if $\Gamma$ is geometrically infinite, we evaluate $\omega_\Gamma$ at the sequence
of points $\{x_n\}$ given by Lemma 5.2,
$$
\omega_\Gamma (x_n) \leq C e^{-\sqrt{\frac {4(\lambda_0 - 2\delta)}{4+\delta}}d_H(
x_n, S_\Gamma)} \text{vol}(B_1(x_n))^{-\frac 12}
\int_{S_\Gamma}\text{vol}(B_1(y))^{-\frac 12}d\sigma(y).
\tag 5.14
$$
Therefore, by Lemma 5.1 and Lemma 5.2, $\omega_\Gamma (x_n) \rightarrow 0$ as $n\rightarrow
\infty$. In light of the formula (4.8), this implies that $n$-dimensional Lebesgue
measure of $L(\Gamma)$ has to be positive, which is a contradiction. Thus the proof of this
theorem is finished.
\enddemo

\vskip .2in

\noindent \S 6. {\bf Conformal finiteness} \vskip .1in

This section is concerned with the question: when is a Kleinian
group conformally finite? We would like to give some geometric
criteria for a Kleinian group to be conformally finite. The idea
still is that the hypersurface $G(\Gamma)/\Gamma$ with the metric
induced from the hyperbolic metric is the right geometric
representative for the Kleinian manifold $\Omega(\Gamma)/\Gamma$.
We first observe:

\proclaim{Theorem 6.1} Suppose that $\Gamma$ is a nonelementary
Kleinian group and $G(\Gamma)$ is the Lipschitz graph constructed
in Section 1. Then $\Gamma$ is conformally finite if and only if
the volume of the hypersurface $G(\Gamma)/\Gamma$ in the
hyperbolic manifold $B^{n+1}/\Gamma$ is finite.
\endproclaim

\noindent {\bf Remark 6.1.} \quad This condition has an analogue
in another formulation of geometric finiteness: the thick part of
the convex core be compact. The latter is equivalent to say that
some neighborhood of the convex core for the hyperbolic manifold
$B^{n+1}/\Gamma$ has finite volume (cf.  [6] [15]). \vskip 0.1in

\demo{Proof of Theorem 6.1} Let us begin with the thick-thin decomposition of hyperbolic
manifolds
$B^{n+1}/\Gamma$ with respect to a small number $\epsilon$ which is smaller than the
Margulis constant in the same dimension.
The hypersurface $G(\Gamma)/\Gamma$ is
also decomposed into thick part $W_\epsilon$ and thin part $S_\epsilon$. Clearly at each point
in the thick part $W_\epsilon$ there is the hyperbolic geodesic ball $B_{\frac 12 \epsilon}$
where $G(\Gamma)\bigcap B_{\frac 12 \epsilon}$ belongs to some fundamental domain for
$\Gamma$ on the graph $G(\Gamma)$. Because $G(\Gamma)$ is Lipschitz graph over the unit
sphere (or any sphere with the same center), the volume of $G(\Gamma)\bigcap B_{\frac 12
\epsilon}$ under the metric induced from the hyperbolic metric on $B^{n+1}$ is bounded from below
by some constant only depending on $\epsilon$. Therefore, if $G(\Gamma)/\Gamma$ has a finite
volume with the metric, then the thick part $W_\epsilon$ has to be compact. Now we may conclude
that the number of the noncompact connected components has to be finite. Because the finite
boundary of each noncompact component, which is the connecting region
of the end to the thick part,
has a size again bounded from below (depending on $\epsilon$). Notice that, each of those noncompact
thin ends corresponds to a maximum parabolic subgroup $P_i$ whose fixed point is $p_i$. Then,
we find that, for some fundamental domain for $\Gamma$ in its domain of discontinuity
$\Omega(\Gamma)$, there are only finite number of limit points $p_i$ on its boundary. Moreover,
those parabolic fixed points therefore have to be bounded, i.e. have to be so-called cusped
limit points. This means precisely that the Kleinian manifold $\Omega(\Gamma)/\Gamma$ is a disjoint
union of a compact part and a finite number of standard conformal cusp ends. So $\Gamma$ is conformally
finite.

On the other hand if $\Gamma$ is conformally finite, it follows from Theorem 3.5 that
the hypersurface $G(\Gamma)/\Gamma$ has finite volume. So the proof is completed.
\enddemo

As a consequence we have the following criterion to tell when a Kleinian group is
conformally finite.

\proclaim{Theorem 6.2} Suppose that $\Gamma$ is a nonelementary
Kleinian group. Then $\Gamma$ is conformally finite if and only if
the Kleinian manifold $\Omega(\Gamma)/\Gamma$ possesses a
conformal metric such that

(1) its volume is finite;

(2) $|R| + |\nabla R| \leq C$ and $Ric \geq -C$.
\endproclaim

\demo{Proof} First of all, if $\Gamma$ is conformally finite, it
is clear the Kleinian manifold possesses a metric satisfying (1)
and (2), in light of the discussion in Section 3. The converse
part of this theorem is a consequence of above Theorem 6.1 and
Theorem 2.12, Chapter VI in [16]. We remark that, a stronger
assumption that $M$ has bounded curvature is listed for the above
result, but it is clear from the proof (e.g. applying method of
gradient estimate), that assumptions as (2) given here are
sufficient for the conclusion. More precisely, on the domain of
discontinuity $\Omega(\Gamma)$ the metric is complete and
satisfies (2). The Harnack estimate in [16] shows that, denoting
the metric by $e^{2u}g_o$ where $g_0$ is the standard metric on
the sphere,
$$
u(x) \geq C \frac 1{d(x)} \ \ \text
{for all $x\in \Omega (\Gamma)$ }.
$$
where $d(x) = \text{dist}(x, \partial \Omega (\Gamma))$.
Then for the metric induced from the hyperbolic metric on the hypersurface
$G(\Gamma)/\Gamma$, its volume is controlled by the volume of the
Kleinian manifold with the given metric, therefore is finite. In light of
the above Theorem 6.1, the proof is finished.
\enddemo

\noindent {\bf Remark 6.2.} \quad Conditions (1) and (2) in
Theorem 6.2 can be replaced by a simpler one: finite volume and
all curvature and their derivatives bounded.

\noindent
\vskip .3in
\centerline{\bf References}

\bigskip
\roster

\vskip .1in
\item"{[1]}" Ahlfors, Mobius Transformations in Several Dimensions. Lecture Notes, School of
Math., Univ. of Minnesota, Minneapolis, 1981.

\vskip .1in
\item"{[2]}" Ahlfors, Finitely generated Kleinian groups.
American Journal of Math., vol. 86 (1964), 413-429.

\vskip .1in
\item"{[3]}" L. Bers, On Ahlfors' finiteness
theorem.  Amer. J. Math 89 (1967), 1078-1082.

\vskip 0.1in
\item"{[4]}" L. Bers, Inequalities for finitely
generated Kleinian groups.  J. Analyse Math. 18 (1967), 23-41.

\vskip .1in
\item"{[5]}" C. Bishop and P. Jones, Hausdorff dimension and Kleinian groups.
Acta Math., 179 (1997), 1-39.

\vskip .1in
\item"{[6]}" Bowditch, Geometric finiteness for hyperbolic groups.  J. Funct. Anal. 113
(1993), 245-317.

\vskip .1in
\item"{[7]}" S.-Y. A. Chang,  J. Qing and P. Yang, On the Chern-
Gauss-Bonnet Integral for conformal metrics on $\Bbb R^4$. Duke
Journal of Math, Vol. 103, No. 3, (2000), 523-544.

\vskip .1in
\item"{[8]}" S.-Y. A. Chang,  J. Qing and P. Yang,
Compactification of a class of conformally flat 4-manifolds.
 Inventiones Math. 142 no. 1 (2000),  65--93.

\vskip .1in
\item"{[9]}" Davies, Heat Kernels and Spectral Theory. Cambridge Tracts in Math., 92. Cambridge
Univ. Press, Cambridge-New York, 1989.

\vskip .1in
\item"{[10]}" P. J\"{a}rvi and M. Vuorinen, Uniformly perfect sets and quasiregular
mappings. J. London Math. Soc. (2) 54 (1996), 515529.

\vskip .1in
\item"{[11]}" M. Kapovich, On the absence of Sullivan's cusp
finiteness theorem in higher dimensions.  Algebra and analysis
(Irkutsk, 1989), 77--89, Amer. Math. Soc. Transl. Ser. 2, 163.

\vskip .1in
\item"{[12]}" M. Kapovich and L. Potyagailo, On the absence
of finiteness theorems of Ahlfors and Sullivan for Kleinian groups in
higher dimensions. (Russian) Sibirsk. Mat. Zh. 32 (1991), no. 2,61--73,
212 translation in Siberian Math. J. 32 (1991), no. 2, 227--237.

\vskip .1in
\item"{[13]}" P. Nicholls, The Ergodic Theory of Discrete Groups. London Math. Soc. Lecture Note
Ser., 143. Cambridge Univ. Press, Cambridge, 1989.

\vskip 0.1in
\item"{[14]}" Hiromi Ohtake, On Ahlfors' weak
finiteness theorem.  J. Math Kyoto Univ. 24-4 (1984), 725-740.

\vskip .1in
\item"{[15]}" J. Ratcliff, Foundations of hyperbolic manifolds.
1994 Springer Verlag.

\vskip .1in
\item"{[16]}" R.Schoen and S.T. Yau, Lectures on Differential
Geometry. International Press, 1994.

\vskip .1in
\item"{[17]}" Sullivan, Related aspects of positivity in Riemannian
geometry. J. Diff. Geom., 25(1987), 327-351.

\endroster
\enddocument